\documentclass[12pt]{article}
\usepackage[dvips]{graphicx}
\DeclareGraphicsExtensions{ps,eps}
\makeatletter
\def\section{\@startsection{section}{1}{\z@}{-1.5ex plus -.5ex
minus -2.ex}{1ex plus .2ex}{\large\bf}}

\input{amssym.def}
\input{amssym}

\def\@thmcounterstep{}
\long\def\@makecaption#1#2{\vskip 10pt \setbox\@tempboxa\hbox{#1.#2}
\ifdim \wd\@tempboxa >\hsize
#1.#2\par 
\else99
\hbox to\hsize{\hfil\box\@tempboxa\hfil}
\fi}
\def\ps@headings{
\def\@oddhead{\footnotesize\rm\hfill\runninghead\hfill}
\def\@evenhead{\@oddhead}
\def\@oddfoot{\rm\hfill\thepage\hfill}\def\@evenfoot{\@oddfoot}}
\newtheorem{Theorem}{Theorem}[section]

\makeatother
\title
{
REMARKS
ON THE MINIMIZING GEODESIC PROBLEM
IN INVISCID INCOMPRESSIBLE FLUID MECHANICS
}
\def\runninghead{\quad 
The minimizing geodesic problem in fluid mechanics}
\author{
{\em Yann Brenier}\thanks{CNRS, Universit\'e de Nice Sophia-Antipolis (FR 2800 W. D\"oblin),
Institut Universitaire de France}}
\date{} 
\begin{document}
\pagestyle{headings}
\flushbottom
\maketitle
\vspace{-10pt}

\section{Abstract}

We consider $L^2$ minimizing geodesics along the
group of volume preserving maps $SDiff(D)$ of a given 3-dimensional domain $D$. The corresponding curves describe the motion of an ideal incompressible fluid inside $D$ and are (formally) solutions of the Euler equations. It is known that there is a unique possible pressure gradient for these curves whenever their end points are fixed. In addition, this pressure field has a limited but unconditional (internal) regularity.
The present paper completes these results by showing: 1) the uniqueness property can be viewed as an infinite dimensional phenomenon (related to the possibility of relaxing the corresponding minimization problem by convex optimization) , which is false for finite dimensional configuration spaces such as $O(3)$ for the motion of rigid bodies;
2) the unconditional partial regularity is necessarily limited.
\subsubsection*{}
Key words: calculus of variations, geodesics, fluid mechanics, global analysis
\\
MSC: 35Q35 (49J45, 49N60)
\section{Introduction}

There are very few problems in mathematical fluid mechanics for which global existence and uniqueness results are known without any restriction on the data. (For instance, the Leray
global existence result for 3D Navier-Stokes equations lacks uniqueness, the Yudovich theorem 
for 2D Euler equations requires bounded vorticity, the Glimm-Bressan theory for 1D gas dynamics is only valid for small initial
data, etc...). One of them is the problem of minimizing geodesics ('shortest paths')
on the group of 3D volume
preserving diffeomorphisms with $L^2$ metric, which,
following Arnold's geometric interpretation \cite{AK} of the Euler equations,
is a secondary way to find solutions of these equations without solving an initial value problem.
Let us be more specific.
\\
Given a smooth enough bounded domain $D$ in $R^3$, we denote by $SDiff(D)$ the group
of all volume and orientation preserving diffeomorphisms of $D$.  
This group is a natural configuration space for the motion of an incompressible fluid
moving inside $D$ which is mathematically described by a time dependent curve
$t\rightarrow g_t\in SDiff(D)$.
The minimizing geodesic (or shortest path) problem can be defined as follows:
given two such diffeomorphisms
$g_0$, $g_1$, find a curve $t\in [0,1]\rightarrow g_t\in SDiff(D)$
that achieves the geodesic distance between $g_0$ and $g_1$:
\begin{equation}
\label{distance}
\Delta(g_0,g_1)=\inf\;\sqrt{\int_{0}^{1}||\frac{dg_t}{dt}||_{L^2}^2 dt}
\end{equation}
where $||\cdot||_{L^2}$ denotes the norm in $L^2(D,R^3)$. The formal optimality equation
reads
\begin{equation}
\label{euler}
\frac{d^2 g_t}{dt^2}
\circ g_t^{-1}+\nabla p_t=0,
\end{equation}
where
$p_t$ is a time dependent scalar field defined on $D$ (called the 'pressure field')
which balances the incompressibility constraint as a Lagrange multiplier.
Equation (\ref{euler}) is precisely the Euler equation introduced in 1755 \cite{Eu}
to describe the motion of an inviscid incompressible fluid moving inside $D$ without
any external force. Let us empasize that the minimizing geodesic problem is
different from the more conventional Cauchy problem, for which the initial
'velocity'  $\frac{dg_0}{dt}$ is prescribed together with $g_0$ (which is traditionnaly
normalized to be the identity map) and, of course, there is no prescribed
endpoint $g_1$ at time $t=1$. The minimizing geodesic (or shortest path) problem
has been solved at the local level (together with the Cauchy problem)
by Ebin and Marsden \cite{EM}: if $g_1
\circ g_0^{-1}-I$ is sufficiently small in Sobolev norm
$H^{s}$, $s>5/2$, there is a unique minimizing geodesic.
In sharp contrast, a striking result of Shnirelman \cite{Sh1} shows that the minimizing geodesic problem
may have no solution at all in the large. In the same paper, it is proven that, as a metric
space, the completion of $SDiff(D)$ for the geodesic distance (\ref{distance}) is the semi-group $VPM(D)$ of
all volume-preserving maps $g$ of $D$ in the measure theoretic sense: $g$ is Borel and
$$
\int_D f(g(x))dx=\int_D f(x)dx,
$$
holds true for all continuous functions $f$.  Notice that, in particular, any map
$g(x)=(h(x_1),x_2,x_3)$ , for which $h$ is a Lebesgue measure-preserving map
of the unit interval $[0,1]$, belongs to the completion of $SDiff(D)$ when $D=[0,1]^3$!
However, even in the completed configuration space $VPM(D)$, the minimizing
geodesic problem may have no solution, as shown by the author in \cite{Br1} with
the example:
$$
D=[0,1]^3,\;\;\;g_0(x),\;\;\;g_1(x)=(h(x_1),x_2,x_3),\;\;\;h(x_1)=1-x_1.
$$
More positively,  it was proven in \cite{Br2,Br3}
that there is a unique pressure field attached to the minimizing geodesic problem.
With the help of a crucial density result obtained by Shnirelman \cite{Sh2}, this result can be 
stated in the following way:

\begin{Theorem}
\label{global}
Let  $g_0$ and $g_1$ given in $VPM(D)$ where $D=[0,1]^3$.
We say that $g^\epsilon_t\in SDiff(D)$ is an approximate minimizing geodesic
between $g_0$ and $g_1$ if
\begin{equation}
\label{approximate geodesic}
\int_{0}^{1}||\frac{dg^\epsilon_t}{dt}||_{L^2}^2 dt \rightarrow \Delta(g_0,g_1)^2,
\;\;\;||g^\epsilon_0-g_0||_{L^2}^2\rightarrow 0,
\;\;\;||g^\epsilon_1-g_1||_{L^2}^2\rightarrow 0,
\end{equation}
as $\epsilon\downarrow 0$.
Then, there is a unique pressure gradient field $\nabla p_t$ such that, for all approximate minimizing geodesic,
$$\frac{d^2 g^\epsilon_t}{dt^2}\circ (g_t^\epsilon)^{-1}+\nabla p_t\rightarrow 0$$
holds true in the sense of distributions in the interior of $[0,1]\times D$.
\end{Theorem}

The main idea of the proof is that
the minimizing geodesic problem, which apparently is a minimization problem
lacking both compactness and convexity, can be relaxed as
a CONVEX minimization problem, in an appropriate generalized framework,
without relaxation gap, thanks to Shnirelman's density result.
This is strongly related to the fact that
the completed configuration space $VPM(D)$ is 
a dense subspace (as shown in \cite{Ne,BG}) of a weakly compact convex set, namely $DS(D)$, the set of all doubly stochastic measures on $D\times D$ (i.e. all nonnegative
measures on $D\times D$
having the Lebesgue measure as projection on both copies of $D$), with respect to
the weak convergence of measures, through the
embedding $g\in VPM(D)\rightarrow \mu_g\in DS(D)$, where
$$
\int_{D^2}f(x,y)d\mu_g(x,y)=\int_D f(x,g(x))dx,\;\;\;\forall f\in C(D\times D).
$$
The uniqueness of $p$ follows, quite easily, from a duality argument due to this 'hidden' convex structure, as shown in \cite{Br2}. 
Nevertheless, geometrically speaking, the uniqueness of $p$ is quite surprising, 
since, between two given points, minimizing
geodesics are not necessarily unique (as can be easily checked) whenever they exist
and there are no a priori reasons that the corresponding acceleration fields $-\nabla p$
should be identical. It is unlikely that such a property could be proven using classical differential geometric tools.
\\
Theorem \ref{global} can be completed in various manners, as in \cite{Br3} (existence of
generalized minimizing geodesics solving a generalized version of the Euler
equations, partial regularity of $p$, etc.) or, recently,
by Ambrosio and Figalli  in \cite{AF1,AF2}, Bernot, Figalli, Santambrogio  \cite{BFS}. 
In particular, $\nabla p$ is shown to be
a locally bounded measure in the interior of $[0,1]\times D$ in \cite{Br3}.
The time integrability of $p$ has been improved since by Ambrosio and Figalli
\cite{AF2}: $p$ belongs to the space $L^2_{loc}(]0,1[,BV_{loc}(D^\circ))$.
At this point, two natural questions can be asked:
\\
Q1: Is the uniqueness of the pressure a specific property of minimizing geodesics on the
infinite dimensional group of volume preserving diffeomorphisms,
which is not true for similar finite dimensional configuration spaces,
such as the special orthogonal group $SO(3)$ for the motion of a rigid body?
\\
Q2: Is the regularity of the pressure field only partial?
\\
In both cases, the present paper provides a positive answer.

\section{About the uniqueness of the pressure field}

In this section, we give some evidence that the uniqueness of the pressure field for the minimizing
geodesic problem on $SDiff(D)$ (or its completion $VPM(D)$) is a genuine
infinite-dimensional phenomenon. For this purpose, we consider the finite dimensional situation
where the special orthogonal group $SO(3)$ substitutes for $SDiff(D)$ and rigid motions of solid bodies substitute for incompressible
inviscid fluid motions. (See  \cite{Ar}.)
As a matter of fact, as shown below, rigid motions
can be interpreted as particular solutions of the Euler equations, in the special case
when the fluid domain $D$
is an ellipsoid $D=KB_1$, where $B_1$ is the unit ball in $R^3$ and $K$  a symmetric
positive matrix.
To see that, we first notice that a volume preserving map $g:D=KB_1\rightarrow D$ 
is a linear map if and only if:
\begin{equation}
\label{rigid}
g(x)=KUK^{-1}x,\;\;\;x\in D=KB_1,
\end{equation}
for some matrix  $U$ in the orthogonal group  $O(3)$. Indeed, by definition of
$D=KB_1$, $x\rightarrow K^{-1}g(Kx)$ must be a linear map preserving the unit ball $B_1$,
i.e.  an element of $O(3)$.
The condition $U\in O(3)$ can be expressed in a variational way by requiring
\begin{equation}
\label{constraint}
Tr(U^*MU)=Tr(M)
\end{equation}
for all  $3\times 3$ real symmetric matrices $M$,
where $Tr$  and $*$ respectively denote
the trace and the transposition operator on $3\times 3$ real matrices.
\\
Let us now look at  geodesics $(g_t)$
on $SDiff(D)$ with the additional constraint that $g_t$ must be linear at each time
$t$ ('rigid motions'), i.e. of form (\ref{rigid}): $g_t=KU_tK^{-1}$ for some curve $U_t$
valued in $O(3)$. 
We have
$$
\int_D |\frac{dg_t(x)}{dt}|^2 dx
=\int_D  |K\frac{dU_t}{dt}K^{-1}x|^2 dx
$$
$$
=r_0\int_{B_1}  |K\frac{dU_t}{dt}x|^2 dx
=r _1\; Tr (\frac{dU_t^*}{dt}K^2\frac{dU_t}{dt})
$$
where $r_0,r_1>0$ are normalization factors.
Thus, 
encoding the condition $U\in O(3)$,
the geodesics we are looking for correspond to  saddle points 
$(U_t,M_t;\;t\in [0,1])$ of the Lagrangian
\begin{equation}
\label{lagrangian}
\int_{0}^{1} \{Tr (\frac{dU_t^*}{dt}K^2\frac{dU_t}{dt})-Tr(U_t^*M_tU_t)+tr(M_t)\}dt.
\end{equation}
Here $M_t$ is constrained to be symmetric and there is no more constraint on $U_t$.
The optimality conditions are straightforward:
\begin{equation}
\label{top}
K^2\frac{d^2U_t}{dt^2}+M_tU_t=0,
\end{equation}
with $U_t$ orthogonal and $M_t$ symmetric.
At this point, we have exactly recovered the usual equations for motions of a rigid body
in classical Mechanics (see \cite{Ar}). They describe geodesic curves on $SO(3)$
for the metric generated by the 'inertia' matrix $K^2$ of the body.
\\
The point of our discussion
is that rigid motions not only can be embedded in the framework of fluid Mechanics but,
more strikingly, are just special solutions of the Euler equations.
Indeed, from (\ref{top}), we recover a solution of the Euler equations just by
setting
\begin{equation}
\label{dictionary}
g_t(x)=KU_tK^{-1}x,\;\;\;p_t(x)=\frac{1}{2}K^{-1}M_tK^{-1}x\cdot x,\;\;\;x\in D=KB_1.
\end{equation}
We conclude that geodesics on $SO(3)$ are just special geodesics on
$SDiff(D)$ for a special choice of the domain $D$ (namely the ellipsoid
$D=KB_1$ associated to the inertia matrix $K^2$). 
\\
Now, we can discuss the issue of $minimizing$, not just plain, geodesics.
The key point is that, for a solution of the rigid body motion (\ref{top}) and the
associated solution (\ref{dictionary}) of the Euler equations, which corresponds
to a geodesic curve on both $SO(3)$ and $SDiff(D)$, it is NOT equivalent to
be a $minimizing$ geodesic on $SO(3)$ and $SDiff(D)$.  This makes sense, since
rigid motions are more restrained than fluid motion. So it is conceivable that
$U_t$ is a minimizing geodesic on $SO(3)$, meanwhile
the corresponding $g_t=KU_tK^{-1}x$ is a plain, non minimizing,  geodesic on $SDiff(D)$.  
This suggests that the uniqueness of the pressure field for fluid motions has no equivalent property in the case of rigid motions. As a matter of fact, this follows
from the following (stronger) result:
\begin{Theorem}
\label{SO3}
Let $(U_t,M_t)$, $(V_t,N_t)$ two distinct solutions of the rigid motion equations
for $t\in [0,1]$ such that $U_0=V_0=I$ (where $I$ is the identity matrix) and
$U_1=V_1$. Assume that $M_t=N_t$ for all $t\in [0,1]$. Then $U_t$ and $V_t$
are exceptional in the sense that they
must be rigid rotations with constant angular speed around one of the inertia
axis of the body (i.e. an eigenvector of the inertia matrix $K^2$).
\end{Theorem}

In other words, generically, two distinct geodesics meeting at two different points
must have different accelerations. This is, a fortiori, also true for $minimizing$
geodesics.  The group $SO(3)$ is a smooth closed bounded set in the finite dimensional 
Hilbert space of real $3\times 3$ matrices and has a finite geodesic diameter.
Thus, every geodesic curve can be minimizing only on finite time intervals,
and, therefore,
they are plenty of distinct minimizing geodesics connecting a same pair of points.

\subsection*{Proof of theorem \ref{SO3}}
Let us consider two distinct minimizing geodesics on $SO(3)$,
$(U_t,M_t)$, $(V_t,N_t)$,
with same endpoints $U_0=I$, $U_1$ at $t=0$ and $t=1$.
It is classical to introduce two fields of skew symmetric real matrices,  $B_t$ and $C_t$,
such that 
$$
\frac{dU_t}{dt}=B_t U_t,\;\;\;
\frac{dV_t}{dt}=C_t U_t
$$
(this is like introducing Eulerian coordinates in fluid Mechanics).
Then, the motion equation (\ref{top}) can be written
\begin{equation}
\label{top2}
\frac{dB_t}{dt}+B_t^2=-K^{-2}M_t,
\end{equation}
and, similarly
$$
\frac{dC_t}{dt}+C_t^2=-K^{-2}N_t.
$$
Let us assume that $M_t$ and $N_t$ coincide, while $B_t$ and $C_t$ are
distinct.
Since $B_t^2$ and $C_t^2$ are symmetric meanwhile $B_t$ and $C_t$ are
skew symmetric, it follows that
$C_t=B_t+L$ for some constant skew symmetric matrix $L$ different from zero,
meanwhile  $C_t^2=B_t^2$. This implies $B_tL+LB_t+L^2=0$. and, therefore,
$$
\frac{dB_t}{dt}L+L\frac{dB_t}{dt}=0.
$$
 Up to a change of orthonormal frame, we may assume, without loss of
 generality,
 $Lx=(x_2,-x_1,0)\beta$ for some constant $\beta\ne 0$.
We deduce, by direct calculation, that $\frac{dB_t}{dt}=0$.
So $B_t=B_0$, $C_t=B_0+L$. Using again $B_tL+LB_t+L^2=0$.
we get no other solution than $B_0x=(-x_2,x_1,0)\beta/2$ and $C_0=-B_0$.
This means that $U_t=\exp(B_0 t)$ and $V_t=\exp(-B_0 t)$ are just rotations at constant angular speed
$+\beta/2$ and $-\beta/2$, along the axis $(0,0,1)$. Notice that, since the time interval
for which the geodesics are minimizing has been fixed to be $[0,1]$,
the only possibility is $\beta^2=4\pi^2$. Going back to  (\ref{top2}), we further deduce
$M_t x=M_0x=-K^2B_0^2 x=K(x_1,x_2,0)\pi^2$ where $M_t$ is supposed to be symmetric.
This implies that the axis of rotation $(0,0,1)$ must be an eigenvector for the inertia
matrix $K^2$ (i.e. an inertia axis for the rigid body), 
which is clearly an exceptional situation, as soon as $K$ is a generic
symmetric positive matrix (which corresponds to a generic rigid ellipsoid).

\section{ Limited regularity of the pressure field}

In this section, we provide an explicit, self-similar, solution of the minimizing
geodesics with limited regularity.

\begin{Theorem}
\label{regularity}
Let  $L>0$, $D=[-L,L]\times[0,1]^2$, $L\ge 1$.
\\
Then,  for $0\le t\le 1$, $x\in D$,
\begin{equation}
\label{solution1}
g_t(x)=t^{2/3}(2\sqrt{x_1t^{-2/3}}-1,x_2,x_3),\;\;\;\;
0<x_1<t^{2/3},
\end{equation}
\begin{equation}
\label{solution2}
g_t(x)=t^{2/3}(1-2\sqrt{-x_1t^{-2/3}},x_2,x_3),\;\;-t^{2/3}<x_1<0,
\end{equation}
\begin{equation}
\label{solution3}
g_t(x)=x,\;\;\;|x_1|>t^{2/3},
\end{equation}
define a (generalized) minimizing geodesic, with a pressure field of limited regularity
\begin{equation}
\label{solution}
p_t(x)=p_t(x_1)=-\frac{1}{9t^2}(t^{4/3}-x_1^2)_+.
\end{equation}
\end{Theorem}

\subsection*{Remarks}
i) The family $(g_t,\;0\le t\le 1)$ is not valued in $SDiff(D)$ but in its completion
$VPM(D)$. So, our example does not prevent a better regularity of the pressure
field in the case of smooth data $g_0$, $g_1$, valued in $SDiff(D)$. 
However, it does rule out 
unlimited  internal regularity, in the style of classical elliptic PDE theory, independently on
the boundary data.
This example also shows, in our opinion, that the regularity to be expected for the pressure field
is semi-concavity, or, at least, measure-valued second order space derivatives. Let us recall that, so
far, we only know that $\nabla p$ is  a locally bounded measure \cite{Br5} (or more precisely $p$ belongs to
$L^2_{loc}(]0,1[,BV_{loc}(D^{\circ}))$ \cite{AF2}). So there should be one order of
differentiability in space to be gained in the future.
\\
ii) This solution can be interpreted as a 'hydrostatic vortex sheet' as explained in 
section \ref{appendix hydrostatic} and, under that form, coincides with a self-similar  'relaxed solution' of the Euler equations already introduced by Duchon and Robert in \cite{DR}.

\subsection*{Proof}
The proof is based on two statements:
\\
i) For each $t$, $g_t$ is a volume preserving map (which is not obvious at first glance).
\\
ii) For almost every fixed $x\in D$, $\xi_t=g_t(x)$ minimizes
\begin{equation}
\label{action}
 \int_{0}^{1}
\{\frac{1}{2}|\frac{d\xi_t}{dt}|^2-(p_t)(\xi_t)\}dt,
\end{equation}
among all curve $\xi_t$ valued in $D$ such that $\xi_0=g_0(x)$,
$\xi_1=g_1(x)$.
These two conditions guarantee that, indeed, $g_t$ is a minimizing geodesic,
following a standard argument (see \cite{Br3,AF1}, for instance).
\\
Let us first check the second statement. If $|x_1|>1$, $g_t(x)=x$
and the statement is trivial. If $|x_1|\le 1$, and $x_1\ne 0$, we see that
$\xi_t=g_t(x)$ is continuously differentiable in $t\in [0,1]$, with piecewise continuous second order derivative (with a jump at $t=|x_1|^{3/2}$), and satisfies
\begin{equation}
\label{acceleration}
\frac{d^2\xi_t}{dt^2}=-(\nabla p_t)(\xi_t)
\end{equation}
for almost every $t$, which guarantees that $\xi$ is already a critical point of the action defined by
(\ref{action}). Let us now prove that $\xi$ is also a global minimizer of the action (\ref{action}).
For this purpose, we
compute the second variation $SV$ of 
action (\ref{action}) 
for a perturbation
$\xi_t+\zeta_t$ with $\zeta_0=\zeta_1=0$, and want to show that $SV\ge 0$. We find
$$
SV=
\frac{1}{2}\{ \int_{0}^{1}|\frac{d\zeta_t}{dt}|^2-(\partial^2 p_t)((\xi_t)_1)(\zeta_t)_1^2\}dt.
$$
By  definition (\ref{solution}) of $p_t$,  
$$
\partial^2 p_t(s)=
\frac{2}{9t^2}
[
1\{|s|<t^{2/3}\}-t^{2/3}\delta(s-t^{2/3})-t^{2/3}\delta(s+t^{2/3})
]
$$
Thus
$$
SV\ge 
\frac{1}{2} \int_{0}^{1}|\frac{d\zeta_t}{dt}|^2 dt-
\frac{1}{9t^2} \int_{|(\xi_t)_1|<t^{2/3}}(\zeta_t)_1^2 dt,
$$
Let us recall  the classical Hardy inequality (in one space dimension):
\begin{equation}
\label{hardy}
 \int_{0}^{+\infty}
\{(\frac{d\eta_t}{dt})^2-\frac{\eta_t^2}{4t^2}dt\}\ge 0,
\end{equation}
for all smooth real function  $\eta_t$  such that $\eta_0=0$.
We deduce that $SV\ge 0$ and conclude that $\xi$, indeed, is  a minimizer
of (\ref{action}) as its end points are fixed at $t=0$ and $t=1$.
\\
Finally, let us prove that $g_t$ is volume preserving. Due to the self-similarity of $g_t$,
it is enough to check, in the case $t=1$, that, for every real continuous  function $f$, 
$$
\int_{-1}^{0} f(1-2\sqrt{-x})dx+\int_{0}^{1} f(2\sqrt{x}-1)dx=
\int_{-1}^{1} f(x)dx
$$
which also follows from elementary calculations.
\\
Finally, let us mention that the solution discussed in this section has a fluid mechanic
interpretation and can be derived from the hydrostatic limit of the Euler equations.
See all details in section \ref{appendix hydrostatic}.

\section{Appendix on the hydrostatic limit of the Euler equations}
\label{appendix hydrostatic}

As explained in \cite{Br6}, 
the minimizing geodesic problem is strongly linked to the hydrostatic limit of the Euler
equations, which reads, on the domain $D=[-L,L]\times [0,1]^2$,
\begin{equation}
\label{euler hydrostatic}
D_t v_1+\partial_1 p=0,\;\;\;
D_t v_2+\partial_2 p=0,\;\;\;
D_t=\partial_t+v\cdot\nabla,
\end{equation}
\begin{equation}
\label{euler hydrostatic2}
\partial_3 p=0,\;\;\; \nabla\cdot v=0,\;\;\;v//\partial D.
\end{equation}
These equations are formally obtained  by ignoring
the vertical acceleration term $D_t v_3$ in the classical Euler equations
(see \cite{Li,Br4,Gr,Br5} for some rigorous results).
Notice that, given a sufficiently smooth solution  $(v,p)$ of these hydrostatic equations,
we may introduce the corresponding flow $X(t,x)$, defined by
\begin{equation}
\label{flow}
\partial_t X(t,x)=v(t,X(t,x)),\;\;\; X(t=0,x)=x,
\end{equation}
which provides a time dependent family of  volume-preserving maps $X(t,\cdot)$
of $D$, since $v$ is divergence-free and parallel to the boundary.
We are going to construct an explicit solution to these equations and, as an output,
the solution (\ref{solution1},\ref{solution2},\ref{solution3},\ref{solution}) used in
Theorem \ref{regularity}.

\subsection*{Construction of an explicit solution}

We first define a divergence-free velocity field (with trivial second component)
by setting: $v_2=0$, $v_1=\partial_3\psi$, $v_3=-\partial_1\psi$,
 where $\psi$ is the 
'stream-function' defined by three different formulae in the domain $D$,
depending on the location. We  first set:
\begin{equation}
\label{stream1}
\psi(t,x_1,x_3)=\frac{t^{-1/3}x_3(\xi-1)}{3},\;\;\;\;
0<x_3<\frac{1+\xi}{2},\;\;\;|\xi|<1,
\end{equation}
where we use the rescaled coordinate $\xi=x_1t^{-2/3}$. Next:
\begin{equation}
\label{stream2}
\psi(t,x_1,x_3)=\frac{t^{-1/3}(x_3-1)(\xi+1)}{3},\;\;\;\;
\frac{1+\xi}{2}<x_3<1,\;\;\;|\xi|<1,
\end{equation}
and, finally:
\begin{equation}
\label{stream3}
\psi(t,x_1,x_3)=0,\;\;\;|\xi|>1.
\end{equation}
Notice that the stream-function is  continuous at the interfaces $\xi=1$, $\xi=-1$
and $x_3=(1+\xi)/2$ and we can easily recover the velocity field $v$ by differentiating $\psi$:
\begin{equation}
\label{velocity1}
v_1=\partial_3\psi
=\frac{t^{-1/3}(\xi-1)}{3},\;\;\;v_2=0,\;\;\;v_3=-\partial_1\psi
=-\frac{t^{-1}x_3}{3}
\end{equation}
whenever $0<x_3<\frac{1+\xi}{2},\;\;|\xi|<1$, 
\begin{equation}
\label{velocity2}
v_1
=\frac{t^{-1/3}(\xi+1)}{3},\;\;\;v_2=0,\;\;\;v_3=-\frac{t^{-1}(x_3-1)}{3}
\end{equation}
whenever $\frac{1+\xi}{2}<x_3<1,\;\;|\xi|<1$ ,
and $v=0$ whenever $|\xi|>1$.
\\
This velocity field is piecewise smooth, with a strong singularity at $t=0$
and also at the interfaces $\xi=1$, $\xi=-1$, $x_3=(1+\xi)/2$ for each $t>0$.
The interface $x_3=(1+\xi)/2$ can be interpreted as a vortex sheet initially located
verically above $x_1=0$. This velocity field was advocated by Duchon and Robert
\cite{DR} as an example of 'relaxed solution' of the Euler equations.
\\
Notice the apparent separation of space variables: $v_1$ depends only on $x_1$ (through
$\xi=x_1t^{-2/3}$), $v_3$ depends only on $x_3$ while $v_2=0$.  Strictly speaking,
this is not true, since formulae  (\ref{velocity1},\ref{velocity2}) depend on the sign
of $x_3-\frac{1+\xi}{2}$. However this is good enough to provide a very simple structure
to  the corresponding flow $X(t,x)$ (defined by (\ref{flow}):
$$
X(t,x)=(X_1(t,x_1),x_2,X_3(t,x_1,x_3)),
$$
where the first component depends only on the first space variable (which is not
the case of the last component).
More precisely,
by integration of  (\ref{velocity1}),
we get the following explicit formula:
\begin{equation}
\label{solution1is}
X_1(t,x_1)=t^{2/3}(2\sqrt{x_1t^{-2/3}}-1),\;\;\;\;
0<x_1<t^{2/3},
\end{equation}
\begin{equation}
\label{solution2bis}
X_1(t,x_1)=t^{2/3}(1-2\sqrt{-x_1t^{-2/3}}),\;\;-t^{2/3}<x_1<0,
\end{equation}
\begin{equation}
\label{solution3bis}
X_1(t,x_1)=x_1,\;\;\;|x_1|>t^{2/3}.
\end{equation}
Since $v$ is divergence free and parallel to $\partial D$, 
$X(t,\cdot)$ is a volume preserving map
of $D$. As a consequence; $x_1\rightarrow X_1(t,x_1)$ must be a Lebesgue measure
map of the interval $[-L,L]$. Indeed, for each continuous function $f(x)$,
$$
\int_{D} f(x)dx=\int_{D} f(X(t,x))dx=
\int_{D} f(X_1(t,x_1),x_2,X_3(t,x_1,x_3))dx_1dx_2dx_3
$$
and, in particular, when $f=f(x_1)$:
$$
\int_{-L}^L f(x_1)dx_1=\int_{-L}^L f(X_1(t,x_1))dx_1.
$$
\\
From the definition of $v_1$, we also deduce
\begin{equation}
\label{simple hydrostatic}
\partial_t+\partial(\frac{v_1^2}{2})+\partial_1 p=0,
\end{equation}
(in the sense of distribution, with no spurious singular measure), where
\begin{equation}
\label{simple pressure}
p(t,x_1)=-\frac{1}{9t^2}(t^{4/3}-x_1^2)_+.
\end{equation}
Thus $(v_1,v_2,v_3,p)$ solves the hydrostatic equations (in distribution form).
\\
Also notice that (\ref{simple hydrostatic}) just means:
$$\partial^2_{tt}X_1(t,x_1))=-(\partial_1 p)(t,X_1(t,x_1)).$$
Finally, by setting 
$$
g_t(x)=(X_1(t,x_1),x_2,x_3),\;\;\;p_t(x)=p(t,x_1),
$$
we recover the solution discussed in Theorem \ref{regularity}.

\subsubsection*{Acknowledgments}
The author acknowledges the support of ANR contract OTARIE
ANR-07-BLAN-0235. Part of his research was done during 
the a stay at  UBC, Vancouver.
(Thematic program and PDE summer school,
workshop on Regularity problems in hydrodynamics, organized
by V. Sverak and T.-P. Tsai, PIMS-UBC August 10-14, 2009.)

\end{document}